\documentclass[12pt]{article}
\usepackage{amssymb, amsmath}

\newtheorem{theorem}{Theorem}[section]

\newtheorem{definition}{Definition}[section]

\newtheorem{cor}{Corollary}[section]

\newcommand{\n}{\nonumber}

\newcommand{\si}{\sigma_R (|x|)}

\newcommand{\s}{\sigma}

\newcommand{\bb}{\begin{equation}}
\newcommand{\ee}{\end{equation}}
\newcommand{\bq}{\begin{eqnarray}}
\newcommand{\eq}{\end{eqnarray}}
\newcommand{\bqn}{\begin{eqnarray*}}
\newcommand{\eqn}{\end{eqnarray*}}

\begin{document}
\title{ Liouville type of theorems with weights for the  Navier-Stokes equations
and the Euler equations}
\author{ Dongho Chae\\
 Department of Mathematics\\
  Sungkyunkwan
University\\
 Suwon 440-746, Korea\\
 e-mail: {\it chae@skku.edu }}
 \date{}
\maketitle
\begin{abstract}
We study Liouville type of theorems for the Navier-Stokes and the
Euler equations
 on $\Bbb R^N$, $N\geq 2$. Specifically, we prove that if
a weak solution $(v,p)$ satisfies $|v|^2 +|p| \in L^1 (0,T; L^1(\Bbb
R^N, w_1(x)dx))$ and $\int_{\Bbb R^N} p(x,t)w_2 (x)dx \geq0$ for
some weight functions $w_1(x)$ and $w_2 (x)$, then the solution is
trivial, namely $v=0$ almost everywhere on $\Bbb R^N \times (0, T)$.
Similar results hold for the MHD Equations on $\Bbb R^N$, $N\geq3$.
\end{abstract}
\section{Introduction}
 \setcounter{equation}{0}
We are concerned  on  the Navier-Stokes equations(the Euler
equations for $\nu=0$) on  $\Bbb R^N$, $N\in \Bbb N, N\geq 2$.
\[
\mathrm{ (NS)_\nu}
 \left\{ \aligned
 &\frac{\partial v}{\partial t} +(v\cdot \nabla )v =-\nabla p +\nu \Delta
 v +f
 \quad (x,t)\in \Bbb R^N\times (0, \infty) \\
 & \textrm{div }\, v =0 , \quad (x,t)\in \Bbb R^N\times (0,
 \infty)\\
  &v(x,0)=v_0 (x), \quad x\in \Bbb R^N
  \endaligned
  \right.
  \]
  where $v(x,t)=(v^1 (x,t), \cdots, v^N (x,t))$ is the velocity, $p=p(x,t)$ is the
  pressure,  $f=(f_1(x,t), \cdots, f^N (x,t))$ is the external force, and $\nu \geq 0$ is the viscosity.
Given $a,b \in \Bbb R^N$, we denote by $a\otimes b $  the $N\times
N$ matrix with $(a\otimes b)_{ij}=a_ib_j$.
  For two $N\times N$ matrices $A$ and $B$ we denote $A:B=\sum_{i,j=1}^N A_{ij} B_{ij}$.
 Given $m\in \Bbb N\cup \{0\}, q\in [1, \infty]$, we introduce
 $$W^{m,q}_\s (\Bbb R^N):= \left\{ v\in [W^{m,q} (\Bbb R^N)]^N,\,\,\mathrm{ div}\, v=0 \right\},
 $$
 where $W^{m,q} (\Bbb R^N)$ is the standard Sobolev space on $\Bbb R^N,$ and
 the derivatives in  the operation of div $(\cdot)$ are in the sense of
 distribution. In particular, $H^m_\s (\Bbb R^N):=W^{m,2}_\s (\Bbb
 R^N)$ and $L^q_\s (\Bbb R^N):=W^{0, q}_\s (\Bbb R^N)$.
 Similarly, given $q\in [1, \infty]$, we use $L^q _{ loc, \s } (\Bbb R^N)$ to denote the class of solenoidal
 vector fields, which belongs to $[L^q _{loc} (\Bbb R^N)]^N$.
  In $\Bbb R^N$ we define  weak
 solutions of the Navier-Stokes(Euler)
equations as follows.
  \begin{definition} We say that a pair $(v,p)\in L^2 (0, T; L^2_{loc,\s} (\Bbb R^N))\times L^1
  (0, T; L^1_{loc}(\Bbb R^N ))$ is a weak solution of $(NS)_\nu$
  on $ \Bbb R^N \times(0, T )$ with $f\in L^1(0, T; L^1_{loc, \s}
  (\Bbb R^N))$
  if
  \bq\label{11}
  \lefteqn{-\int_{0}^{T} \int_{\Bbb R^N} v(x,t) \cdot \phi(x)\xi ' (t)dxdt -\int_{0}^{T}\int_{\Bbb R^N} v(x,t)
  \otimes v(x,t):\nabla \phi (x) \xi (t) dxdt}\hspace{.in} \n \\
  && = \int_{0} ^{T}\int_{\Bbb R^N} p(x,t)\mathrm{ div }\, \phi (x)
  \xi (t)dxdt +\nu \int_{0}^{T} \int_{\Bbb R^N} v(x,t)\cdot \Delta \phi (x) \xi (t)dxdt\n \\
  &&\quad +\int_{0}^{T} \int_{\Bbb R^N} f(x,t)\cdot \phi (x) \xi (t)dxdt\n \\
  \eq
   for all $\xi \in C^1 _0 (0, T)$ and $\phi =[C_0 ^\infty (\Bbb R^N )]^N
   $.
  \end{definition}
In \cite{cha1} it is proved that if a weak solution  $(v, p)$ of the
Euler or Navier-Stokes equations satisfy
 \bb\label{12}
v\in L^2 (0, T; L^2 _\s(\Bbb R^N))\quad \mbox{and}\quad p\in L^1 (0,
T; \mathcal{H}^q (\Bbb R^N))
 \ee
  for some $q\in (0, 1]$, where $\mathcal{H}^q (\Bbb
R^N)$ denotes the Hardy space on $\Bbb R^N$,  then $v(x,t)=0$ almost
everywhere on $\Bbb R^N \times (0, T)$. Furthermore, if $p\in L^1(0,
T; L^1 (\Bbb R^N))$, then there happens the {\em equipartition of
energy} over each component(\cite{cha1}),
$$\int_{\Bbb R^N} v^j(x,t)v^k dx=-\delta_{jk}\int_{\Bbb R^N} p(x,t)dx.
$$
The main purpose of this paper is to further develop the idea
initiated in \cite{cha1} to  obtain  substantially extended
Liouville type of theorems with suitable weight functions for the
associated integrations for the Navier-Stokes equations, the Euler
equations  on $\Bbb R^N$, $N\geq 2$, and the (both viscous and
invicid) MHD equations on $\Bbb R^N$, $N\geq 3$. To the author's
knowledge there exist a previous study on the Liouville type of
theorems in for the 3D Navier-Stokes equations  with {\em
axisymmetry} for $\nu>0$ (\cite{koc}), which is in completely
different fashion from that of \cite{cha1} and from those studied in
this paper.
 In the case of the Euler equations and the
 MHD equations, in particular, there exists
 no previous Liouvillle type of results available in the literature.  Our first main theorem is the
following.
 \begin{theorem} Let $w\in L^1_{loc} ([0, \infty))$ be given,
 which is positive almost everywhere on $[0, \infty)$. Suppose  $(v,p)$ is a weak solution to $(NS)_\nu$
 with $f\in L^1(0, T; L^1_{loc, \s} (\Bbb R^N))$ and $\nu \geq 0$  on $\Bbb R^N \times (0, T)$ such that
\bq\label{13}
 \lefteqn{\int_0 ^T \int_{\Bbb R^N} (|v(x,t)|^2 +|p(x,t)|)\times}\hspace{.0in}\n \\
 && \times\left[w(|x|) +\frac{1}{|x|}\int_0 ^{|x|} w(s)ds +\frac{1}{|x|^2}
 \int_0 ^{|x|} \int_0 ^r
 w(s)dsdr \right] dxdt
 <\infty,\n \\
 \eq
 and
 \bb\label{14}
 \int_{\Bbb R^N} p(x,t)\left[w(|x|) +\frac{N-1}{|x|}\int_0 ^{|x|} w(s)ds
 \right]dx \geq 0\quad \mbox{for $t\in (0, T)$}.
 \ee
 Then, $v(x,t)=0$ almost everywhere on $\Bbb R^N \times (0, T)$.
\end{theorem}
{\it Remark 1.1 } If we choose $w(s)\equiv1$ on $[0, \infty)$, then
we recover  Liouville part of results of Theorem 1.1 (i) in
\cite{cha1}. \\
\ \\
\noindent{\it Remark 1.2 } Let us set $w^*(r):=\sup_{0\leq s\leq r}
w(s)$. Then, since
$$ w(r) +\frac{1}{r}\int_0 ^{r} w(s)ds
+\frac{1}{r^2}
 \int_0 ^{r} \int_0 ^s
 w(\rho)d\rho ds\leq\frac{5}{2} w^* (r),$$
 we can replace the condition (\ref{13}) by a stronger one,
 \bb\label{15}
  \int_0 ^T \int_{\Bbb R^N} (|v(x,t)|^2 +|p(x,t)|)w^* (|x|) dxdt
 <\infty
  \ee
  to get our conclusion of the theorem from (\ref{14}).\\
  \ \\
The following is a consequence of the above theorem, which we state
as a separate theorem.
\begin{theorem} Let  $(v,p)$ be a weak solution to $(NS)_\nu$
 with $f\in L^1(0, T; L^1_{loc, \s} (\Bbb R^N))$ and $\nu \geq 0$
 on $\Bbb R^N \times (0, T)$ such that either
\bb\label{16}
 \int_0 ^T \int_{\Bbb R^N}\frac{(|v(x,t)|^2 +|p(x,t) |)}{1+|x|}dxdt <\infty,
 \ee
 or
\bq\label{17a}
&&\mbox{$p(x,t)\to 0$ as $|x|\to \infty$ for almost every $t\in (0, T)$, and }\\
\label{17b}
  &&v\in L^2(0, T; L^q (\Bbb R^N))\qquad \mbox{for some
$q$ with }\quad
 2<q<\frac{2N}{N-1}.
 \eq
Suppose there exists $w\in L^1 (0, \infty)$ such that
 \bb\label{17c}
 0< w(r)\leq \frac{C}{1+r} \quad \mbox{for some $C>0$}
 \ee
 for almost every $r\in [0, \infty)$, and
 \bb\label{18}
 \int_{\Bbb R^N} p(x,t)\left[w(|x|) +\frac{N-1}{|x|}\int_0 ^{|x|} w(s)ds
 \right]dx \geq 0\quad \mbox{for almost every $t\in (0, T)$}.
 \ee
 Then, $v(x,t)=0$ almost everywhere on $\Bbb R^N \times (0, T)$.
\end{theorem}
{\it Remark 1.3 } The main novelty of the above theorem, compared to
Theorem 1.1, is that the integrability conditions (\ref{16}) and
(\ref{17a}) do not involve restriction on the weight function
$w(r)$. Moreover, we do not need any integrability condition on
pressure $p(x,t)$ in (\ref{17a}). The price to pay for theses
relaxations is that we need to select weight functions from a
smaller class than of Theorem 1.1.\\
\ \\
Since $H^1(\Bbb R^N)\hookrightarrow L^{\frac{2N}{N-2}} (\Bbb R^N)$
by the Sobolev embedding,  and $L^2 (\Bbb R^N) \cap
L^{\frac{2N}{N-2}} (\Bbb R^N) \subset L^q
 (\Bbb R^N)$ for $2<q<\frac{2N}{N-2}$ by the standard $L^q(\Bbb R^N)$
 interpolation inequality, we easily find that Leray's weak
 solution(\cite{ler}) to $(NS)_\nu, \nu >0$ satisfies
\bb\label{emb}
 v\in L^\infty (0, T; L^2_\s(\Bbb R^N)) \cap L^2 (0, T; H^1_\s
(\Bbb R^N))\subset L^2 (0, T; L^q _\s (\Bbb R^N))
 \ee
  \quad
for all $q\in (2, \frac{2N}{N-2})$. Hence, as an  immediate
corollary of Theorem 1.2 we obtain:
\begin{cor}
Let $v \in L^\infty (0, T; L^2_\s(\Bbb R^N)) \cap L^2 (0, T; H^1_\s
(\Bbb R^N))$ be  Leray's weak solution to $(NS)_\nu$ with $f\in
L^2(0, T; L^2_{\s} (\Bbb R^N))$ and $\nu >0$. Suppose the pressure
$p(x,t)$ satisfies (\ref{17a}) and (\ref{18}) for a function $w(r)$
satisfying the conditions of Theorem 1.2. Then, $v(x,t)=0$ almost
everywhere on $\Bbb R^N \times (0, T)$.
\end{cor}

The proofs of Theorem 1.1 and Theorem 1.2 are given in the next
section. Further generalized theorems extending them to the MHD
equations are stated and proved in Section 3.

\section{Proof of the Main Theorems}
 \setcounter{equation}{0}
\noindent{\bf Proof of Theorem 1.1 }
 Let us
consider a radial cut-off function $\sigma\in C_0 ^\infty(\Bbb R^N)$
such that
 \bb\label{20}
   \sigma(|x|)=\left\{ \aligned
                  &1 \quad\mbox{if $|x|<1$}\\
                     &0 \quad\mbox{if $|x|>2$},
                      \endaligned \right.
 \ee
and $0\leq \sigma  (x)\leq 1$ for $1<|x|<2$. We set
 \bb\label{21}
W(\rho):=\int_0 ^{\rho} \int_0 ^s w(r)drds.
 \ee
Then, for each $R >0$, we define
 \bb\label{22}
\varphi_R (x)=W(|x|)\s \left(\frac{|x|}{R}\right)=W(|x|)\s_R
(|x|)\in C_0 ^\infty (\Bbb R^N).
 \ee
 Let $\xi \in C^1 _0 (0, T)$, and we   choose the vector test function
 $\phi$ in (\ref{11}) as
 \bb\label{23}
  \phi= \nabla \varphi_R (x).
 \ee
 Then, after routine computations (\ref{11}) becomes
 \bqn
\lefteqn{0=\int_0 ^T\int_{\Bbb R^N} \left[W^{\prime\prime} (|x|)
\frac{(v\cdot x)^2}{|x|^2} + W^{\prime} (|x|)
\left(\frac{|v|^2}{|x|}
-\frac{(v\cdot x)^2}{|x|^3}\right) \right] \si \xi(t)  \,dx dt}\hspace{.5in}\n\\
&&\quad+\int_0 ^T\int_{\Bbb R^N}  W' (|x|) \s'
\left(\frac{|x|}{R}\right)  \frac{(v\cdot x)^2}{R|x|^2}\xi(t) \,dx
dt\n \\
    &&\quad+ \int_0 ^T\int_{\Bbb R^N}
\frac{1}{R}\left( \frac{|v|^2}{|x|} -\frac{(v\cdot x)^2}{|x|^3}
\right) \s'\left(\frac{|x|}{R}\right)W(|x|)\xi(t) \,dx dt
 \eqn
  \bq\label{24}
&&\quad+\int_0 ^T\int_{\Bbb R^N} \frac{(v\cdot x)^2}{ R^2|x|^2}
\s^{\prime\prime} \left(\frac{|x|}{R}\right) W(|x|)\xi(t) \,dx
dt\n \\
 &&\quad+ \int_0 ^T\int_{\Bbb R^N}p(x,t)\left[ W^{\prime\prime}
 (|x|) +(N-1)\frac{W' (|x|)}{|x|}\right]\sigma_R (|x|)
 \xi(t) \, dxdt\n \\
 &&\quad+ \frac{2}{R}\int_0 ^T\int_{\Bbb R^N}p(x,t) W' (|x|)
 \s' \left(\frac{|x|}{R}\right)\xi(t) \, dxdt\n \\
&&\quad+ \frac{N-1}{R}\int_0 ^T\int_{\Bbb R^N}p(x,t)\frac{1}{|x|}\s'
\left(\frac{|x|}{R}\right) W(|x|)\xi(t) \, dxdt\n \\
&&\quad+ \int_0 ^T\int_{\Bbb R^N}p(x,t)\frac{1}{R^2}
\s^{\prime\prime} \left(\frac{|x|}{R}\right)W(|x|)\xi(t) \,
dxdt \n \\
 &&:=I_1+\cdots +I_8
\eq
 Note that the term involving derivative with respect to time, the viscosity term and the forcing term
 in (\ref{11})
 vanish altogether, since
 $$
 \int_0 ^T\int_{\Bbb R^N} v(x,t)\cdot \nabla\varphi_R (x) \xi'(t) dxdt=0,
 $$
$$
 \int_0 ^T\int_{\Bbb R^N} v(x,t)\cdot \nabla (\Delta \varphi_R (x)) \xi(t)
 dxdt=0
 $$
 for  $v\in L^2(0, T; L^2_{loc, \s} (\Bbb R^N))$
 and
 $$
 \int_0 ^T\int_{\Bbb R^N} f(x,t)\cdot \nabla\varphi_R (x) \xi'(t) dxdt=0,
 $$
 for  $f\in L^1(0, T; L^1_{loc, \s} (\Bbb R^N))$
  by the divergence free condition in the sense of distribution.
 In terms of the function $W$ defined in  (\ref{21}) our condition (\ref{13}))
 can be written as
\bb\label{25}
 \int_0 ^T \int_{\Bbb R^N} (|v(x,t)|^2 +|p(x,t)|) \left[W^{\prime\prime}
 (|x|) +\frac{1}{|x|}W' (|x|) +\frac{1}{|x|^2} W(|x|)\right]
 dxdt<\infty.
 \ee
Since
\bqn
 &&\int_0 ^T \int_{\Bbb R^N} \left|\left[W^{\prime\prime} (|x|)
\frac{(v\cdot x)^2}{|x|^2} + W^{\prime} (|x|)
\left(\frac{|v|^2}{|x|} -\frac{(v\cdot x)^2}{|x|^3}\right)
\right]\right||\xi(t)|dxdt\\
&&\qquad \leq 2 \sup_{0\leq t\leq T}
  |\xi(t)|\int_0 ^T\int_{\Bbb R^N} |v(x,t)|^2\left[ W^{\prime\prime}
(|x|)+\frac{W'(|x|)}{|x|} \right] dxdt <\infty, \eqn
  We can use the dominated convergence theorem to show that
  \bb\label{26}
  I_1 \to \int_0 ^T\int_{\Bbb R^N} \left[W^{\prime\prime} (|x|)
\frac{(v\cdot x)^2}{|x|^2} + W^{\prime} (|x|)
\left(\frac{|v|^2}{|x|} -\frac{(v\cdot x)^2}{|x|^3}\right) \right]
\xi(t) \,dx dt
 \ee
  as $R\to \infty$.
Similarly,
  \bb\label{27}
  I_5\to \int_0 ^T\int_{\Bbb R^N}p(x,t)\left[ W^{\prime\prime}
 (|x|) +(N-1)\frac{W' (|x|)}{|x|}\right]
 \xi(t) \, dxdt
\ee as $R\to \infty$.
 For $I_2$ we estimate
 \bq\label{28}
  |I_2 |&\leq& \int_0 ^T \int_{R<|x|<2R} |v(x,t)|^2\left|\s'
\left(\frac{|x|}{R}\right)\right|
  \frac{W'(|x|)}{|x|} \frac{|x|}{R}|\xi (t)|dxdt\n \\
  &\leq &2 \sup_{1<s<2} |\s'(s)|\sup_{0\leq t\leq T}
  |\xi(t)|\int_0 ^T \int_{R<|x|<2R} |v(x,t)|^2
  \frac{W'(|x|)}{|x|} dxdt\n \\
\to 0
  \eq
 as $R\to \infty$ by the dominated convergence theorem.
Similarly
  \bq\label{29}
   |I_3|&\leq &2 \int_0 ^T\int_{R<|x|<2R}  \frac{|x|}{R}|v(x,t)|^2
\left|\s'\left(\frac{|x|}{R}\right)\right|\frac{W(|x|)}{|x|^2}\xi(t) \,dx dt\n \\
&\leq &4\sup_{1<s<2} |\s'(s)|\sup_{0\leq t\leq T}
  |\xi(t)|\int_0 ^T \int_{R<|x|<2R} |v(x,t)|^2
  \frac{W'(|x|)}{|x|}dxdt
\to 0,\n \\
  \eq
  and
  \bq\label{210}
  |I_4|&\leq&\int_0 ^T\int_{R<|x|<2R}\frac{|x|^2}{R^2} |v(x,t)|^2\left|\s^{\prime\prime}
  \left(\frac{|x|}{R}\right)\right|\frac{W(|x|)}{|x|^2}\xi(t) \,dx dt\n \\
  &\leq&4\sup_{1<s<2} |\s^{\prime\prime}(s)|\sup_{0\leq t\leq T}
  |\xi(t)|\int_0 ^T\int_{R<|x|<2R}|v(x,t)|^2
  \frac{W(|x|)}{|x|^2}\,dx dt\to 0\n \\
  \eq
   as $R\to \infty$. The estimates for $I_6,I_7$ and $I_8$ are
   similar to the above, and we find
   \bq\label{211}
   |I_6|&\leq &2 \int_0 ^T\int_{R<|x|<2R}|p(x,t)| \frac{|x|}{R}\frac{W'
   (|x|)}{|x|}
 \left|\s' \left(\frac{|x|}{R}\right)\right| |\xi(t)| \, dxdt\n \\
 &\leq& 4\sup_{1<s<2} |\s' (s)|\sup_{0\leq t\leq T}
  |\xi(t)|\int_0 ^T\int_{R<|x|<2R}|p(x,t)|\frac{W'
   (|x|)}{|x|}dx dt\to 0,\n \\
  \eq
\bq\label{212}
   |I_7|&\leq & (N-1)\int_0 ^T\int_{R<|x|<2R}|p(x,t)|\frac{|x|}{R}\left|\s'
\left(\frac{|x|}{R}\right)\right| \frac{W(|x|)}{|x|^2}|\xi(t)| \, dxdt\n \\
 &\leq& 2\sup_{1<s<2} |\s' (s)|\sup_{0\leq t\leq T}
  |\xi(t)|\int_0 ^T\int_{R<|x|<2R}|p(x,t)|\frac{W(|x|)}{|x|^2}dx dt\to 0,\n \\
  \eq
  and
  \bq\label{213}
 |I_8|&\leq&\int_0 ^T\int_{\Bbb R^N}|p(x,t)|\frac{|x|^2}{R^2}
\left|\s^{\prime\prime}\left(\frac{|x|}{R}\right)\right|
\frac{W(|x|)}{|x|^2}|\xi(t)| \, dxdt \n\\
  &\leq& 4\sup_{1<s<2}
|\s^{\prime\prime} (s)|\sup_{0\leq t\leq T}
  |\xi(t)|\int_0 ^T\int_{R<|x|<2R}|p(x,t)|\frac{W
   (|x|)}{|x|^2}dx dt\to 0\n \\
  \eq
  as $R\to \infty$ respectively. Thus passing $R\to \infty$ in (\ref{24}),
  we finally obtain
\bq\label{214}
 &&\int_0 ^T\int_{\Bbb R^N} \left[W^{\prime\prime} (|x|)
\frac{(v\cdot x)^2}{|x|^2} + W^{\prime} (|x|)
\left(\frac{|v|^2}{|x|} -\frac{(v\cdot x)^2}{|x|^3}\right) \right]
\xi(t) \,dx dt\n \\
&&\qquad=-\int_0 ^T\int_{\Bbb R^N}p(x,t)\left[ W^{\prime\prime}
 (|x|) +(N-1)\frac{W' (|x|)}{|x|}\right]
 \xi(t) \, dxdt\n \\
 \eq
for all $\xi\in C^1 _0 (0, T)$, which can be written , in terms of
the function $w(r)$,  as
 \bq\label{215}
 &&\int_{\Bbb R^N} \left[w (|x|)
\frac{(v\cdot x)^2}{|x|^2} + \frac{1}{|x|}\int_0 ^{|x|}w (s) ds
\left(|v|^2 -\frac{(v\cdot x)^2}{|x|^2}\right)\right]
\,dx \n \\
&&\qquad=-\int_{\Bbb R^N}p(x,t)\left[ w
 (|x|) +(N-1)\frac{1}{|x|}\int_0 ^{|x|} w (s)ds \right]\,dx\leq 0\n \\
 \eq
for almost every $t\in (0, T)$ by the hypothesis (\ref{14}). Since
$|v|^2 \geq (v\cdot x)^2/|x|^2$, we need to have \bqn
 &&\int_{\Bbb R^N} w (|x|)
\frac{(v\cdot x)^2}{|x|^2}\, dx=\int_{\Bbb R^N} \frac{1}{|x|}\int_0
^{|x|}w (s) ds \left[|v|^2 -\frac{(v\cdot x)^2}{|x|^2}\right]\, dx =
0
 \eqn
 almost every $t (0, T)$.
 By the hypothesis $w(|x|)>0$ and $\frac{1}{|x|}\int_0
^{|x|}w (s) ds>0$  for almost every $x\in  \Bbb R^N $, and
 we should have $v(x,t)=0$ for almost every $(x,t)\in \Bbb R^N\times (0, T)$.
 $\square$\\
 \ \\
 \noindent{\bf Proof of Theorem 1.2 } The conditions $w\in L^1 (0,
 \infty)$ and (\ref{17c})
 imply that there exists a positive constant
 $C=C(\|w\|_{L^1 (0, \infty)})$ such that
 \bb\label{216}
  w(r)+\frac{1}{r}\int_0 ^r w(s)ds + \frac{1}{r^2}\int_0 ^r \int_0
 ^s w(\rho)d\rho ds \leq \frac{C}{1+r}.
 \ee
 Therefore, if (\ref{16}) holds true, then
 \bq\label{217}
 \lefteqn{ \int_0 ^T \int_{\Bbb R^N} (|v(x,t)|^2 +|p(x,t)|)\times}\hspace{.0in}\n \\
 && \times\left[w(|x|) +\frac{1}{|x|}\int_0 ^{|x|} w(s)ds +\frac{1}{|x|^2}
 \int_0 ^{|x|} \int_0 ^r
 w(s)dsdr \right] dxdt\n \\
 &&\leq C\int_0 ^T \int_{\Bbb R^N} \frac{|v(x,t)|^2
 +|p(x,t)|}{1+|x|} dxdt <\infty
\eq
 Next, we suppose (\ref{17a}) holds true. In this case we have the well-known
  pressure-velocity
 relation
 $$p(x,t)=\sum_{j,k=1}^N R_j R_k ( v_jv_k )(x,t)  $$
  with $R_j, j=1,\cdots N,$
 the Riesz transforms in $\Bbb R^N$, and thus the Calderon-Zygmund
 inequality says(\cite{ste})
 $$\|p(t)\|_{L^{\frac{q}{2}}}\leq C_q \|v(t)\|_{L^q}^2 \quad \forall q\in
 (2, \infty)
 $$
for a constant $C_q$. Hence, for $2<q<\frac{2N}{N-1}$, we can
estimate
 \bqn
 \lefteqn{ \int_0 ^T \int_{\Bbb R^N} (|v(x,t)|^2 +|p(x,t)|)\times}\hspace{.0in}\n \\
 && \times\left[w(|x|) +\frac{1}{|x|}\int_0 ^{|x|} w(s)ds +\frac{1}{|x|^2}
 \int_0 ^{|x|} \int_0 ^r
 w(s)dsdr \right] dxdt
 \eqn
\bq\label{218}
 &&\leq C\int_0 ^T \int_{\Bbb R^N} \frac{|v(x,t)|^2
 +|p(x,t)|}{1+|x|} dxdt \n \\
&&\leq C\int_0 ^T (\|v(t)\|_{L^q}^2 +\|p(t)\|_{L^{\frac{q}{2}}}
 )dt \left(\int_{\Bbb R^N} \frac{dx}{(1+|x|)^{\frac{q}{q-2}}}
 \right)^{\frac{q-2}{q}}\n \\
 &&\leq C \int_0 ^T \|v(t)\|_{L^q}^2 dt <\infty.
 \eq
 Hence, for both of the cases where either (\ref{16}) or (\ref{17a})
 holds true we can applyt Theorem 1.1 to conclude that $v(x,t)=0$ for
almost every $(x,t)\in \Bbb R^N \times (0, \infty)$. $\square$\\

\section{The case of the MHD equations}
 \setcounter{equation}{0}

In this section we extend the previous results on the system
$(NS)_\nu$ to the magnetohydrodynamic equations in $\Bbb R^N$,
$N\geq 3$.
\[
\mathrm{ (MHD)_{\mu,\nu} }
 \left\{ \aligned
 &\frac{\partial v}{\partial t} +(v\cdot \nabla )v =
 (b\cdot\nabla)b-\nabla (p +\frac12 |b|^2) +\nu \Delta v +f, \\
 &\frac{\partial b}{\partial t} +(v\cdot \nabla )b =(b \cdot \nabla
 )v+\mu \Delta b +g,\\
 &\quad \textrm{div }\, v =\textrm{div }\, b= 0 ,\\
  &v(x,0)=v_0 (x), \quad b(x,0)=b_0 (x)
  \endaligned
  \right.
  \]
where $v=(v_1, \cdots , v_N )$, $v_j =v_j (x, t)$, $j=1,\cdots,N$,
is the velocity of the flow, $p=p(x,t)$ is the scalar pressure,
$b=(b_1, \cdots , b_N )$, $b_j =b_j (x, t)$, is the magnetic field,
 and $v_0$, $b_0$ are the given initial
velocity and magnetic field,
 satisfying div $v_0 =\mathrm{div}\, b_0= 0$, respectively.
 We may consider $f=(f_1(x,t), \cdots , f_N (x,t))$ and $g=(g_1(x,t), \cdots , g_N
(x,t))$ as external forces for the velocity and to the magnetic
fields, respectively. If we set $b=g=0$, then $(MHD)_{\mu,\nu}$
reduces to $(NS)_\nu$ of the previous sections. Let us
 begin with the definition of the weak solutions of $(MHD)_{\mu,\nu}$.
   \begin{definition} We say the triple of functions $(v,b,p)\in [L^2 (0, T; L^2_{loc, \s} (\Bbb R^N))]^2\times L^1 (0, T; L^1 _{\mathrm{loc}} (\Bbb R^N ))$ is a weak solution of $ (MHD)_{\mu,\nu}$
  on $ \Bbb R^N \times(0, T )$,
  if
  \bq\label{31}
  \lefteqn{-\int_{0}^{T} \int_{\Bbb R^N} v(x,t) \cdot \phi(x)\xi ' (t)dxdt -\int_{0}^{T}\int_{\Bbb R^N} v(x,t)
  \otimes v(x,t):\nabla \phi (x) \xi (t) dxdt} \n \\
  &&=-\int_{0}^{T}\int_{\Bbb R^N}b(x,t)
  \otimes b(x,t):\nabla \phi (x) \xi (t) dxdt
 +\int_{0} ^{T} \int_{\Bbb R^N}p(x,t)\,\mathrm{ div }\, \phi (x)
  \xi (t)dxdt \n \\
  &&\quad+\frac12 \int_{0} ^{T} \int_{\Bbb R^N}|b(x,t)|^2\,\mathrm{ div }\, \phi (x)
  \xi (t)dxdt
  +\nu \int_{0}^{T} \int_{\Bbb R^N} v(x,t)\cdot \Delta \phi (x) \xi (t)dxdt,\n \\
  &&\qquad + \int_{0}^{T} \int_{\Bbb R^N} f(x,t)\cdot\phi (x) \xi (t)dxdt,\n \\
  \eq
  and
  \bq\label{32}
  \lefteqn{-\int_{0}^{T} \int_{\Bbb R^N} b(x,t) \cdot \phi(x)\xi ' (t)dxdt -\int_{0}^{T}\int_{\Bbb R^N} v(x,t)
  \otimes b(x,t):\nabla \phi (x) \xi (t) dxdt} \n \\
  &&=-\int_{0}^{T}\int_{\Bbb R^N}b(x,t)
  \otimes v(x,t):\nabla \phi (x) \xi (t) dxdt+
  \mu \int_{0}^{T} \int_{\Bbb R^N} b(x,t)\cdot \Delta \phi (x) \xi (t)dxdt\n \\
   &&\qquad + \int_{0}^{T} \int_{\Bbb R^N} g(x,t)\cdot\phi (x) \xi (t)dxdt,\n \\
  \eq
   for all $\xi \in C^1 _0 (0, T)$ and $\phi =[C_0 ^2 (\Bbb R^N )]^N
   $.
  \end{definition}
 We have the following theorem.
 \begin{theorem} We fix $\mu, \nu \geq 0$, $N\geq 3$.
 Let $w\in L^1_{loc} ([0, \infty))$ be given, which is positive,
 non-increasing  on $[0, \infty)$.
 Suppose  $(v,b,p)\in [L^2(0, T ; L^2_{loc, \s} (\Bbb R^N))]^2\times L^1(0, T ; L^1_{loc} (\Bbb R^N)$ is
 a weak solution to $(MHD)_{\mu, \nu}$ with $f,g\in L^1(0, T;
 L^1_{\s, loc} (\Bbb R^N))$
on $\Bbb R^N \times (0, T)$  such that
  \bq\label{34}
 \lefteqn{\int_0 ^T \int_{\Bbb R^N} (|v(x,t)|^2 +|b(x,t)|^2+|p(x,t)|)\times}\hspace{.0in}\n \\
 && \times\left[w(|x|) +\frac{1}{|x|}\int_0 ^{|x|} w(s)ds +\frac{1}{|x|^2}
 \int_0 ^{|x|} \int_0 ^r
 w(s)dsdr \right] dxdt
 <\infty,\n \\
 \eq
 and
 \bb\label{35}
 \int_{\Bbb R^N} p(x,t)\left[w(|x|) +\frac{N-1}{|x|}\int_0 ^{|x|} w(s)ds
 \right]dx \geq 0\quad \mbox{for $t\in (0, T)$}.
 \ee
 Then, $b(x,t)=0,$ and $v(x,t)=0$ almost everywhere on $\Bbb R^N \times (0, T)$.
\end{theorem}
{\it Remark 3.1 } Similarly to the case of Euler equations if  we
choose $w(s)\equiv1$ on $[0, \infty)$, then we recover a part of
Liouville type of result in Theorem 3.1 in \cite{cha1}.\\
\ \\
\noindent{\it Remark 3.2 } Similarly to  Remark 1.2 for
$w^*(r):=\sup_{0\leq s\leq r} w(s)$
 we can replace (\ref{34}) by a stronger assumption,
 \bb\label{36}
   \int_0 ^T \int_{\Bbb R^N} (|v(x,t)|^2 +|b(x,t)|^2+|p(x,t)|)w^* (|x|)
dxdt
 <\infty
  \ee
to derive triviality of the solution from (\ref{34}).\\
\ \\
 \noindent{\bf Proof of Theorem 3.1 } The method of proof is
similar to that of Theorem 1.1, and we will be brief, describing
only essential points. Similarly to (\ref{20})-(\ref{23}) we choose
$\xi \in C_0 ^1 (0, T)$ and the vector test function $\phi =\nabla
\varphi_R $, where
 \bb\label{37}
  \varphi_R (x)=W(|x|)\s
\left(\frac{|x|}{R}\right)=W(|x|)\s_R (|x|)
 \ee
with $ W(|x|)=\int_0 ^{|x|} \int_0 ^s w(r)drds, $ and $\s $ is the
cut-off function defined in (\ref{20}).  Then, we obtain from
(\ref{31}) that
  \bq\label{38}
 \lefteqn{0=\int_0 ^T\int_{\Bbb R^N} \left[W^{\prime\prime} (|x|)
\frac{(v\cdot x)^2}{|x|^2} + W^{\prime} (|x|)
\left(\frac{|v|^2}{|x|} -\frac{(v\cdot x)^2}{|x|^3}\right)
\right]\si
\xi(t) \,dx dt}\n \\
&& -\int_0 ^T\int_{\Bbb R^N} \left[W^{\prime\prime} (|x|)
\frac{(b\cdot x)^2}{|x|^2} + W^{\prime} (|x|)
\left(\frac{|b|^2}{|x|} -\frac{(b\cdot x)^2}{|x|^3}\right)
\right]\si
\xi(t) \,dx dt\n \\
&& + \int_0 ^T\int_{\Bbb R^N}(p(x,t)+\frac12 |b|^2)\left[
W^{\prime\prime}
 (|x|) +(N-1)\frac{W' (|x|)}{|x|}\right]\si
 \xi(t) \, dxdt \n \\
 && +o(1),\eq
where $o(1)$ denotes the sum of the terms vanishing as $R\to
\infty$.  Taking $R\to \infty$ in (\ref{38}), and rearranging the
non-vanishing terms, we find that
 \bq\label{39}
 \lefteqn{\int_0 ^T\int_{\Bbb R^N} \left[W^{\prime\prime} (|x|)
\frac{(v\cdot x)^2}{|x|^2} + W^{\prime} (|x|)
\left(\frac{|v|^2}{|x|} -\frac{(v\cdot x)^2}{|x|^3}\right) \right]
\xi(t) \,dx dt}\n \\
&& +\frac12 \int_0 ^T\int_{\Bbb R^N} W^{\prime\prime} (|x|)
|b|^2\xi(t) \,dx dt\n \\
&& +\int_0 ^T\int_{\Bbb R^N} \left[\frac{1}{|x|}  W^{\prime}
(|x|)- W^{\prime\prime} (|x|) \right]\frac{(b\cdot x)^2}{|x|^2}\xi(t) \,dx dt\n \\
&&+\frac{N-3}{2}\int_0 ^T\int_{\Bbb R^N} |b|^2 \frac{ W^{\prime}
(|x|)}{|x|}\xi(t) \,dx dt\n \\
 &&=-  \int_0 ^T\int_{\Bbb R^N}p(x,t)\left[
W^{\prime\prime}
 (|x|) +(N-1)\frac{W' (|x|)}{|x|}\right]
 \xi(t) \, dxdt \n \\
 \eq
for all $\xi \in C^1_0 (0,T)$. Hence
 \bq\label{310}
 \lefteqn{\int_{\Bbb R^N} \left[W^{\prime\prime} (|x|)
\frac{(v\cdot x)^2}{|x|^2} + W^{\prime} (|x|)
\left(\frac{|v|^2}{|x|} -\frac{(v\cdot x)^2}{|x|^3}\right) \right]
 \,dx }\n \\
&& +\frac12 \int_{\Bbb R^N} W^{\prime\prime} (|x|) |b|^2
 \,dx\n \\
&& +\int_{\Bbb R^N} \left[\frac{1}{|x|}  W^{\prime} (|x|)-
W^{\prime\prime} (|x|) \right]\frac{(b\cdot x)^2}{|x|^2} \,dx\n \\
&&\quad+\frac{N-3}{2}\int_{\Bbb R^N} |b|^2 \frac{ W^{\prime}
(|x|)}{|x|}\,dx \n \\
 &&\qquad =-  \int_{\Bbb R^N}p(x,t)\left[
W^{\prime\prime}
 (|x|) +(N-1)\frac{W' (|x|)}{|x|}\right]
 \, dx \n \\
 \eq
for almost every $t\in (0, T)$. Our assumption (\ref{35}) implies
that the right hand side of (\ref{310}) is non-positive. Since each
integral of the left hand side of (\ref{310}) is non-negative for
$N\geq 3$, we need to have that each term of the left hand side of
(\ref{310}) vanishes for almost every $t\in (0, T)$. The
requirement,
$$\int_{\Bbb R^N} \left[W^{\prime\prime} (|x|) \frac{(v\cdot
x)^2}{|x|^2} + W^{\prime} (|x|) \left(\frac{|v|^2}{|x|}
-\frac{(v\cdot x)^2}{|x|^3}\right) \right]
 \,dx=0
 $$
 implies $v(x,t)=0$ for almost every $(x,t)\in \Bbb R^N \times (0, T)$,
 as we in the proof of Theorem 1.1. Since $w(r)$ is non-increasing
 on $[0, \infty)$, we have
 $$
\frac{1}{|x|}  W^{\prime} (|x|)- W^{\prime\prime} (|x|)
=\frac{1}{|x|} \int_0 ^{|x|} w(r)dr -w(|x|)\geq 0,
$$
and therefore
$$
\int_{\Bbb R^N} \left[\frac{1}{|x|}  W^{\prime} (|x|)-
W^{\prime\prime} (|x|) \right]\frac{(b\cdot x)^2}{|x|^2} \,dx=0.
 $$
Hence the condition $w(r) =W^{\prime\prime} (r)>0$ together with the
fact $\int_{\Bbb R^N} W^{\prime\prime} (|x|) |b|^2
 \,dx=0
 $
implies $b(x,t)=0$ for almost every $(x,t)\in \Bbb R^N\times (0,
T)$.
 $\square$\\
\ \\
Similarly to Theorem 1.2, we can establish  the following:
\begin{theorem} Let  $(v,b,p)$ be a weak solution to $(MHD)_{\mu,\nu}$
 with $\mu, \nu \geq 0$ and $f,g\in L^1(0, T; L^1_{loc, \s} (\Bbb R^N))$ on $\Bbb R^N \times (0, T)$, $N\geq 3$,  such that
 either
\bb\label{311}
 \int_0 ^T \int_{\Bbb R^N}\frac{(|v(x,t)|^2 +|b(x,t)|^2+|p(x,t) |)}{1+|x|}dxdt <\infty,
 \ee
 or
\bq\label{312a}
 &&\mbox{$|p(x,t)|\to 0$ as $|x|\to \infty$ for almost every $t\in (0,
 T)$, and} \\
 \label{312b}
 &&|v|+|b|\in L^2(0, T; L^q (\Bbb R^N))\qquad \mbox{for some $q$ with }\quad
 2<q<\frac{2N}{N-1}. \n \\
 \eq
Suppose there exists $w\in L^1 (0, \infty)$, which is positive,
non-increasing on $[0, \infty)$ such that
 \bb\label{313}
  0<w(r)\leq \frac{C}{1+r} \quad \mbox{for some constant $C>0$}
 \ee
almost every $r\in [0, \infty)$, and
 \bb\label{314}
 \int_{\Bbb R^N} p(x,t)\left[w(|x|) +\frac{N-1}{|x|}\int_0 ^{|x|} w(s)ds
 \right]dx \geq 0\quad \mbox{for  almost every $t\in (0, T)$}.
 \ee
 Then, $v(x,t)=b(x,t)=0$ almost everywhere on $\Bbb R^N \times (0, T)$.
\end{theorem}

In the case of $\mu, \nu >0$ a global in time weak solutions
$(v,b)\in  [L^\infty (0, T; L^2_\s(\Bbb R^N)) \cap L^2 (0, T; H^1_\s
(\Bbb R^N))]^2$ are constructed in \cite{ser}. Hence, using the fact
(\ref{emb}) we have the following:
\begin{cor}
Let $(v, b, p)\in [L^\infty (0, T; L^2_\s(\Bbb R^N)) \cap L^2 (0, T;
H^1_\s (\Bbb R^N))]^2 \times L^1 (0, T; L^1_{loc} (\Bbb R^N))$ be a
weak solution to $(MHD)_{\mu,\nu}$ with $f, g \in L^2(0, T; L^2_{\s}
(\Bbb R^N))$ and $\mu, \nu >0$, constructed in \cite{ser}. Suppose
the pressure $p(x,t)$ satisfies (\ref{17a}) and (\ref{18}) for a
function $w(r)$ satisfying the conditions of Theorem 1.2. Then,
$v=b=0$ almost everywhere on $\Bbb R^N \times (0, T)$.
\end{cor}

{\bf Proof of Theorem 3.2 }
 Similarly to the proof of Theorem 1.2 the conditions
$w\in L^1 (0,
 \infty)$ and (\ref{313}) imply that there exists a positive constant
 $C=C(\|w \|_{L^1 (0, \infty)})$ such that
 $$ w(r)+\frac{1}{r}\int_0 ^r w(s)ds + \frac{1}{r^2}\int_0 ^r \int_0
 ^s w(\rho)d\rho ds \leq \frac{C}{1+r}.
 $$
 Therefore, if (\ref{311}) holds true, then
 \bq\label{315}
 \lefteqn{ \int_0 ^T \int_{\Bbb R^N} (|v(x,t)|^2 +|b(x,t)|^2+|p(x,t)|)\times}\hspace{.0in}\n \\
 && \times\left[w(|x|) +\frac{1}{|x|}\int_0 ^{|x|} w(s)ds +\frac{1}{|x|^2}
 \int_0 ^{|x|} \int_0 ^r
 w(s)dsdr \right] dxdt\n \\
 &&\leq C\int_0 ^T \int_{\Bbb R^N} \frac{|v(x,t)|^2+|b(x,t)|^2
 +|p(x,t)|}{1+|x|} dxdt <\infty
\eq
 Next, we suppose (\ref{312a})-(\ref{312b}) holds true.
 In order to handle this case we observe that taking the divergence operation of the first equation
of $(MHD)_{\mu, \nu}$, we obtain
$$ \Delta( p+\frac12 |b|^2)= \sum_{j,k=1}^N \partial_j\partial_k (b_jb_k) -
\sum_{j,k=1}^N \partial_j\partial_k (v_j v_k). $$
  Therefore
$$
 p(x,t)= -\sum_{j,k=1}^N R_jR_j (b_jb_k) (x,t) +\sum_{j,k=1}^N R_jR_k (v_jv_k)(x,t)
 -\frac12 |b(x,t)|^2+h(x,t),
$$
where $R_j =\partial_j (-\Delta)^{-\frac12}$, $j=1,\cdots, N$, is
the Riesz transform, and $h(x,t)$ is a harmonic function on $\Bbb
R^N$. The condition (\ref{312a}) implies that $h(\cdot,t)=0$ for
almost every $t\in (0, T)$. As in the proof of Theorem 1.2, thanks
to the Calderon-Zygmund inequality, we have \bb\label{316}
 \| p(t)\|_{L^{\frac{q}{2}}} \leq C_q ( \|v(t)\|_{L^{q}}^2 +\|b(t)\|_{L^{q}}^2)
 \quad \forall q\in (2, \infty)
 \ee
 for a constant $C_q$.
 Therefore
 \bq\label{317}
 \lefteqn{ \int_0 ^T \int_{\Bbb R^N} (|v(x,t)|^2 +|b(x,t)|^2+|p(x,t)|)\times}\hspace{.0in}\n \\
 && \times\left[w(|x|) +\frac{1}{|x|}\int_0 ^{|x|} w(s)ds +\frac{1}{|x|^2}
 \int_0 ^{|x|} \int_0 ^r
 w(s)dsdr \right] dxdt\n \\
 &&\leq C\int_0 ^T \int_{\Bbb R^N} \frac{|v(x,t)|^2+|b(x,t)|^2
 +|p(x,t)|}{1+|x|} dxdt \n \\
 &&\leq C \int_0 ^T (\|v(t)\|_{L^q}^2 +\|b(t)\|_{L^q}^2+\|p(t)\|_{L^{\frac{q}{2}}}
 )dt \left(\int_{\Bbb R^N} \frac{dx}{(1+|x|)^{\frac{q}{q-2}}}
 \right)^{\frac{q-2}{q}}\n \\
 &&\leq C \int_0 ^T (\|v(t)\|_{L^q}^2+\|b(t)\|_{L^q}^2)dt <\infty
 \eq
 for $2<q<\frac{2N}{N-1}$, where we used (\ref{316}).
Therefore, for both of the cases wether (\ref{311}) or
(\ref{312a})-(\ref{312b}) holds true, we can apply Theorem 3.1 to
conclude that $v(x,t)=b(x,t)=0$ for
almost every $(x,t)\in \Bbb R^N \times (0, \infty)$. $\square$\\
$$\mbox{\bf Acknowledgements} $$
The author would like to thank Professor Antonio C\'ordoba for
stimulating discussions. This work was supported partially by  KRF
Grant(MOEHRD, Basic Research Promotion Fund)

\end{document}